\documentclass[11pt]{amsart}
\usepackage{amsfonts,amssymb,amscd,amsmath,enumerate,verbatim,calc}

\newcommand{\CM}{Cohen-Macaulay}
\newcommand{\IFF}{\text{if and only if}}
\newcommand{\wrt}{with respect to}

\newcommand{\m}{\mathfrak{m} }

\newcommand{\x}{\mathbf{x}}

\newcommand{\xar}{\longrightarrow}

\newcommand{\depth}{\operatorname{depth}}

\theoremstyle{plain}
\newtheorem{theorem}{Theorem}

\newtheorem{question}[theorem]{Question}

\theoremstyle{definition}

\newtheorem{example}[theorem]{Example}

\theoremstyle{remark}

\begin{document}
\title[Invariance]{Invariance of a length \\ associated to a reduction }
\author{Tony ~J.~Puthenpurakal}

\date{\today}
\address{Department of Mathematics, IIT Bombay, Powai, Mumbai 400 076}

\email{tputhen@math.iitb.ac.in}
 \begin{abstract}
Let $(A, \m)$ be a $d$-dimensional \CM \ local ring with infinite residue field and let
$J$ be a minimal reduction of $\m$. We show that $ \lambda(\m^3/J\m^2)$ is independent of $J$.  
\end{abstract}

 \maketitle
The notion of minimal reduction of an ideal in a local Noetherian ring  $(A,\m)$ was introduced by 
Northcott and Rees \cite{NoR}. 
It has significant applications in the theory of Hilbert
function of $\m$-primary ideals in $A$, particularly in the case when $A$ is \CM.
Minimal reductions of an ideal $I$  are highly non-unique. In fact if the residue field of
$A$ is infinite and $I$ is $\m$-primary then any $d$ general linear combinations
of generators of $I$, (here  $d = \dim A$),  gives a  minimal reduction of $I$.   
 
However when $A$ is \CM \  with infinite residue field there are some invariants of the ring
 and  the ideal $I$ which are 
independent of minimal  reductions of  $I$. To state it let us fix  some notation. Let $I$ be an $\m$-primary
ideal of $A$ and let  $J$ be a minimal  reduction of $I$. 
We let $\lambda(N)$ denote the length of an $A$-module $N$ and $\mu(N)$
number of its minimal generators. Let   $G_{I}(A) =  \oplus_{n\geq 0}I^n/I^{n+1}$
be  the associated graded ring of $A$ \wrt \ $I$.  We let $e_{0}^{I}(A)$ denote the multiplicity
of $A$ \wrt \  $I$. 

It is well known that
\begin{enumerate}
\item
$\displaystyle{\lambda(I/J) = e_{0}^{I}(A) - \lambda(A/I) \quad \text{due to Serre, cf.
 \cite[Theorem 4.7.6]{BH} } }$.
\item
$\displaystyle{ \lambda(I^2/JI) =  e_{0}^{I}(A) + (d -1)\lambda(A/I) -  \lambda(I/I^2) \quad 
\text{\cite[Lemma 1]{Va1}} }$.
\end{enumerate}
Furthermore if $\depth G_I(A) \geq d -1$ then $\lambda(I^{n+1}/JI^n)$ is independent of 
$J$ \cite[Corollary 3.9]{Mr2}. Furthermore in this case the coefficients of the Hilbert polynomial
of $A$ \wrt \ $I$ can also be expressed in terms of $\lambda(I^{n+1}/JI^n)$, with $n \geq 0$ 
\cite[ Corollary 2.1]{Huck2}.  So it is of some interest to see how these lengths behave when we 
do not have any 
assumptions on depth $G_I(A)$. In general $ \lambda(I^3/JI^2)$ is not independent of 
minimal reduction (see Example \ref{ex}). Our result shows that $ \lambda(\m^3/J\m^2)$ is
invariant of minimal reduction $J$ of $\m$. 

The fact essentially used in our proof is the following result: let $I$ be an $\m$-primary ideal of $A$ and let $J$ be a reduction of $I$,
 and suppose that $J$ is minimally generated by $x_1,\ldots,x_n$. Then $J$ is a minimal 
reduction of $I$ \IFF \ the elements  $x_1,\ldots,x_n$ are analytically independent in $I$ 
and $n = \dim A $.
 Recall that $x_1,\ldots,x_n$
are \emph{analytically independent in}  $I$ if whenever $f(X_1,\ldots,X_n)$ is a homogeneous polynomial of
degree $m$ in $A[X_1,\ldots,X_n]$ ($m$ arbitrary) such that $f(x_1,\ldots,x_n) \in I^m\m$ then 
all coefficients  of $f$ are in $\m$. 

\begin{theorem}
Let $(A,\m)$ be a \CM \ local ring of dimension $d \geq 1$ with infinite residue field. If
$J$ is a minimal reduction of $\m$ then we have an equality

$\displaystyle{  \lambda(\m^3/J\m^2) = e+ (d-1)\mu(\m) -  \mu(\m^2)  - \binom{d-1}{2}}.$
\end{theorem}
\begin{proof}
When $d = 1$ we have to show that $\mu(\m^2) + \lambda(\m^3/J\m^2) = e$. This is well known cf. 
\cite[Theorem 6.18]{VaSix}. So we assume that $\dim A \geq 2$. We assert that it suffices to construct an exact sequence

\begin{equation}
\label{suff}
0 \xar \left( \frac{A}{\m} \right)^{\binom{d}{2}}\xrightarrow{\psi_d} \left( \frac{\m}{\m^2} \right)^{d}\xrightarrow{\phi_d}  \frac{J\m}{J\m^2}  \xar 0
\end{equation}
Suppose we have the exact sequence as claimed. Then we prove the result as follows:
The exact sequence (\ref{suff}) gives that
\begin{equation}
\label{suff2}
\lambda \left( \frac{J\m}{J\m^2} \right)  = d\mu(\m) - \binom{d}{2}
\end{equation}
We also have 
\begin{align*}
\lambda(J\m/J\m^2) &= \lambda(A/J\m^2) - \lambda(A/J\m) \\
                   &= \lambda(A/\m^3) + \lambda(\m^3/J\m^2) - \lambda(A/\m^2) - \lambda(\m^2/J\m) \\
                   &= \lambda(\m^2/\m^3 ) + \lambda(\m^3/J\m^2) - \lambda(\m^2/J\m) \\
                   &= \mu(\m^2) + \lambda(\m^3/J\m^2) - (e - (1 + \mu(\m) -d )) 
\end{align*}
So by (\ref{suff2}) we get 
\begin{align*}
 \mu(\m^2) + \lambda(\m^3/J\m^2) &= e - (1 + \mu(\m) -d) + d\mu(\m) -  \binom{d}{2} \\
                                  & = (d-1)\mu(\m) + e + d - 1 -  \binom{d}{2}  \\
                          & = (d-1)\mu(\m) + e   -  \binom{d -1}{2}
\end{align*}
So it remains to construct the exact sequence (\ref{suff}).

\emph{The author thanks the referee for indicating a simpler proof than the original.}

Denote the Koszul complex on $x_1 \ldots,x_d$ by $K_A(\x)$.
Consider part of $K_A(\x)$ 
\[
  A^{\binom{d}{2}} \xrightarrow{\Psi_d} A^d  \xrightarrow{\Phi_d} A. 
\]
Then $\rm{Image}~\Psi_d \subseteq {\m}^{\oplus d}$ and 
$\Phi_d({\m}^{\oplus d})\subseteq J{\m}$, which gives a right exact
sequence 
\[
 A^{\binom{d}{2}}\xar  {\m}^{\oplus d} \xar J{\m} \xar 0. 
\]
Now tensoring this by $A/\m$ we get the sequence (1) which is right exact.

Hence we only need to prove that $\psi_d$ is injective. 
Let $\{e_i\mid 1\leq i\leq d\}$ and $\{e_i\wedge e_j \mid i<j\}$ be
the canonical bases of $(A/{\m})^d$
and $(A/{\m})^{\binom{d}{2}}$ respectively. Let $\sum_{i<j} m_{ij}e_i\wedge
e_j\in \ker \psi_d$ then 
$\psi_d(\sum m_{ij}e_i\wedge e_j) = \sum_in_i e_i$, where

$$n_i = m_{1i}x_1 +\cdots + m_{i-1,i}x_{i-1} -(m_{i,i+1}x_{i+1}+\cdots
m_{d,i}x_d).$$ 

By  analyticity we get 
$m_{ij} \in {\m}$ for all $j$ and $i$.  
Therefore we have constructed an exact sequence (\ref{suff}) as claimed and
as noted before this completes the proof.
\end{proof}

Next we give an example to show that if $I \neq \m$ then the assertion of  
Theorem 1. does not hold. This example was constructed by Huckaba \cite[ 3.1]{H3}
to show that reduction number of $I$ \wrt \ to minimal reduction $J$ of $I$, depends on $J$.
\begin{example}
\label{ex}
Let $A = k[[x,y]]$ and $I = (x^7,x^6y,x^2y^5,y^7)A$. The ideals $J_1 = (x^7,y^7)A$ and 
$J_2 = (x^7, x^6y + y^7)$ are minimal reductions of $I$. Note that all the ideals
involved are homogeneous. So to compute lengths $\lambda(I^3/J_1I^2)$ and $\lambda(I^3/J_2I^2)$
we may work in the polynomial ring $k[x,y]$ and use some computer algebra package (we used
'Singular') to get
\[
\lambda(I^3/J_1I^2) = 3, \quad \text{and} \quad  \lambda(I^3/J_2I^2) = 2.
\]
\end{example}
Our theorem also prompts the following:
\begin{question}
Let $(A,\m)$ be a \CM \ local ring of dimension $d \geq 1$ with infinite residue field. If
$J$ is a minimal reduction of $\m$ then is $ \lambda(\m^4/J\m^3)$ independent of $J$?
\end{question}
Even though I personally think that this question does not have a positive answer, I have not been
able to get a counter-example.

\providecommand{\bysame}{\leavevmode\hbox to3em{\hrulefill}\thinspace}
\providecommand{\MR}{\relax\ifhmode\unskip\space\fi MR }
\providecommand{\MRhref}[2]{%
  \href{http://www.ams.org/mathscinet-getitem?mr=#1}{#2}
}
\providecommand{\href}[2]{#2}

\end{document}